\documentclass[11pt]{article}
\usepackage[margin=1in]{geometry}
\usepackage{xcolor}
\usepackage{amsmath,amssymb,amsthm,mathtools}
\usepackage{enumitem}
\usepackage{hyperref}
\usepackage{mathrsfs}
\usepackage{tikz-cd}
\usepackage{amscd}
\usepackage[nameinlink,capitalise,noabbrev]{cleveref}
\IfFileExists{stmaryrd.sty}{%
  \usepackage{stmaryrd}%
}{%
}
\hypersetup{
  colorlinks=true,
  linkcolor=blue!50!black,
  citecolor=blue!50!black,
  urlcolor=blue!50!black
}

\allowdisplaybreaks
\numberwithin{equation}{section}
\usepackage{microtype}
\newtheorem{theorem}{Theorem}[section]
\newtheorem{proposition}[theorem]{Proposition}

\newtheorem{corollary}[theorem]{Corollary}
\newtheorem{question}[theorem]{Question}
\theoremstyle{definition}
\newtheorem{definition}[theorem]{Definition}
\newtheorem{convention}[theorem]{Convention}
\theoremstyle{plain}
\newtheorem{maintheorem}{Theorem}

\newtheorem{maincorollary}[maintheorem]{Corollary}

\crefname{maintheorem}{theorem}{theorems}
\Crefname{maintheorem}{Theorem}{Theorems}
\crefname{maincorollary}{corollary}{corollaries}
\Crefname{maincorollary}{Corollary}{Corollaries}

\newcommand{\Met}{\operatorname{Met}^{\infty}}
\newcommand{\Vol}{\operatorname{Vol}}
\newcommand{\Area}{\operatorname{Area}}
\newcommand{\Mass}{\mathbf M}
\newcommand{\Flat}{\mathcal F}

\newcommand{\Bd}{\mathcal B}
\newcommand{\Half}{\mathcal H}
\newcommand{\hw}{\widetilde\omega}
\newcommand{\Z}{\mathbb Z}

\newcommand{\Q}{\mathbb Q}
\newcommand{\eps}{\varepsilon}
\DeclareMathOperator{\Ric}{Ric}

\title{\textbf{Interlacing of the Volume and Half-Volume Spectra}}
\author{Xiaoxiang Jiao \and Hangyue Zhu}
\date{\today}
\begin{document}
\maketitle
\begin{abstract}
We give a direct proof of the weak interlacing inequalities \[\omega_p\leq\hw_p\leq\omega_{p+1}\] between the volume and half-volume spectra of a closed connected Riemannian manifold $M^{n+1}$, $n+1\geq3$. We also prove that, for a residual set of smooth metrics, the area of every smooth embedded half-volume CMC boundary lies outside the rational span of the areas of smooth closed embedded minimal hypersurfaces. For dimensions $3\le n+1\le 5$, these results yield strict interlacing for generic metrics and a lower bound of order $A^{n+1}$ for the number of distinct half-volume CMC boundaries with area at most $A$.

\end{abstract}

\tableofcontents

\section{Introduction}

Let $M^{n+1}$ be a smooth, closed, connected manifold of dimension $n+1\geq3$. The volume spectrum, also known as $p$-widths, associates to each Riemannian metric $g$ a sequence $\{\omega_p(M,g)\}_{p\geq1}$, obtained from areas of min-max minimal hypersurfaces detected by $p$-sweepouts \cite{Law}. Mazurowski--Zhou introduced the half-volume spectrum $\{\hw_p(M,g)\}_{p\geq1}$ by requiring every member of the sweepout family to enclose volume $\Vol_g(M)/2$ \cite{half}. In the present work, we study the comparison between these two spectra.

%GMT
We first recall the definitions involved and fix notation. Write $\mathfrak h=\frac12\Vol_g(M)$, and let $\mathcal C(M)$ denote the space of Caccioppoli sets in $M$, equipped with the $L^1$ topology. Let $\Bd(M;\Flat;\Z_2)=\partial\mathcal C(M)$ be the space of mod-$2$ boundaries equipped with the $\mathcal F$-topology, and let $\Half(M;\Z_2)=\{\partial\Omega:\Omega\in\mathcal C(M),\Vol_g(\Omega)=\mathfrak h\}$. We use $\Mass$ to denote mass of mod-$2$ currents. 

%p-sw
An admissible $p$-sweepout is a continuous map $\Phi:X\to\Bd(M;\Flat;\Z_2)$, where $X$ be a finite cubical subcomplex, such that \begin{enumerate}[label=(\roman*)]
    \item $\Phi^*(\lambda^p)\neq0$. Here $\lambda\in H^1(\Bd(M;\mathbb Z_2);\mathbb Z_2)$ denotes the fundamental cohomology class given by the Almgren isomorphism; see \cite{ALMGREN1962257,Law}.
    \item $\Phi$ has no concentration of mass, i.e. $\lim_{r\downarrow0}\sup_{x\in X}\sup_{q\in M}\Mass\bigl(\Phi(x)\llcorner B_g(q,r)\bigr)=0.$
\end{enumerate}
Let $\mathcal P_p(M,g)$ be the collection of these sweepouts, and let $\mathcal Q_p(M,g)$ consist of those taking values in $\Half(M;\Z_2)$. The $p$-th volume and half-volume spectra are
\[
\begin{aligned}
    \omega_p(M,g)&=
    \inf_{\Phi\in\mathcal P_p(M,g)}\sup_{x\in\operatorname{dom}\Phi}\Mass(\Phi(x)),\\
    \hw_p(M,g)&=\inf_{\Phi\in\mathcal Q_p(M,g)}\sup_{x\in\operatorname{dom}\Phi}\Mass(\Phi(x)).
\end{aligned}
\]
We follow the conventions of \cite{Law,half}. Note that $\omega_p\leq\hw_p$ from the definitions. 

In \cite{half}, Mazurowski and Zhou proved the weak interlacing for phase transition half-volume spectrum, i.e. $c(p)\leq\tilde{c}(p)\leq c(p+1)$. They also suggested the following analogue and proposed an approach using the cup product.

\begin{question}[Mazurowski--Zhou]\label{que:interlacing}
Does
\[
    \omega_p(M,g)\leq\hw_p(M,g)\leq\omega_{p+1}(M,g)
\]
hold for every smooth Riemannian metric $g$ on $M$ and every $p\geq1$?
\end{question}

\paragraph{Main results. }Throughout, we use $\Met(M)$ to denote the space of smooth Riemannian metrics with the $C^\infty$ topology, and CMC hypersurface stands for constant mean curvature hypersurface.

\begin{maintheorem}[Theorem \ref{thm:mainprecise}]\label{thm:A}
Let $M^{n+1}$ be a smooth, closed, connected manifold with $n+1\geq3$. The following hold for every $p\geq1$:
\begin{enumerate}[label=(\roman*)]
    \item For every smooth Riemannian metric $g$,
    \[
        \omega_p(M,g)\leq\hw_p(M,g)\leq\omega_{p+1}(M,g).
    \]
    \item If $3\leq n+1\leq5$, there exists a residual set $\mathcal G'\subset\Met(M)$ such that for every
    $g\in\mathcal G'$,
    \[
        \omega_p(M,g)<\hw_p(M,g)<\omega_{p+1}(M,g).
    \]
\end{enumerate}
\end{maintheorem}
The strict inequalities in (ii) use the following generic separation result between CMC and minimal areas. A half-volume CMC boundary is a smooth embedded boundary $\Sigma=\partial\Omega$ with $\Vol_g(\Omega)=\mathfrak h$ and has constant mean curvature. 

\begin{maintheorem}[Theorem \ref{thm:nonres}]\label{thm:B}
Let $M^{n+1}$ be a smooth, closed, connected manifold with $n+1\geq3$. There exists a residual set $\mathcal R\subset\Met(M)$ such that, for every $g\in\mathcal R$ and every smooth embedded half-volume CMC boundary $\Sigma$,
\begin{equation}
    \Area_g(\Sigma)\notin
    \operatorname{span}_{\Q}
    \bigl\{\Area_g(\Gamma):\Gamma\text{ is a smooth closed embedded minimal hypersurface}\bigr\}.
\end{equation}
\end{maintheorem}

In dimensions three to five, Theorem \ref{thm:B} and the existing min-max realization theorems in\cite{Law,MZCMC} yield the strict inequalities in Theorem \ref{thm:A} (ii).

Mazurowski and Zhou already proved generic infinite existence of half-volume CMC hypersurfaces \cite[Corollary~1.2]{MZCMC}. Strict interlacing and the Weyl law \cite[Theorem~1.1]{Law} give the following quantitative bound.

\begin{maincorollary}[Corollary \ref{cor:Cprecise}]\label{cor:C}
For $3\leq n+1\leq5$, $g\in\mathcal G'$, let $\mathcal N_g(A)$ denote the number of geometrically distinct smooth embedded half-volume CMC boundaries with nonzero mean curvature and area at most $A$. Then
\[
    \liminf_{A\to\infty}\frac{\mathcal N_g(A)}{A^{n+1}}
    \geq\frac{1}{a_n^{n+1}\Vol_g(M)^n},
\]
where $a_n$ is the constant in \cite[Theorem~1.1]{Law}. 
\end{maincorollary}

\paragraph{Related works. }Existence, especially multi-existence of CMC hypersurfaces has been an important question in min-max theory. Zhou--Zhu developed the min-max theory for CMC hypersurfaces and proved the existence of a smooth, closed, almost embedded CMC hypersurface in closed ambient manifolds of dimensions three to seven \cite{CMCI}. Dey proved the quantitative multi-existence result of small $c$-CMC hypersurfaces \cite{DeyCMC} while Pacard--Xu gave the quantitive multi-existence result for large $c$-CMC hyperspheres \cite{MR2481045}. Mazurowski--Zhou proved the existence of infinitely many geometrically distinct half-volume CMC hypersurfaces for generic metrics and for metrics of positive Ricci curvature in dimensions three to five \cite{MZCMC}, and this result could be broaden to dimension six via Chen--Hong--Li's result \cite{R6}. For generic metrics in dimensions three to seven, they also established an alternative result: for each $c>0$, either there exists infinitely many closed $c$-CMC hypersurfaces, or there exists infinitely many half-volume CMC hypersurfaces with mean curvature less than $c$ \cite{MZPisa}.

An equivalence between the Almgren--Pitts and phase-transition half-volume spectra is stated in Talipov's recent preprint \cite{Talipov}. Combining this equivalence with the comparison $c(p)=\omega_p$ given by Dey \cite{Dey222} and the inequality $\tilde c(p)\leq c(p+1)$ given by Mazurowski--Zhou \cite[Proposition~24]{half}, would also give the same result with Theorem \ref{thm:A}(i). We prove this upper bound through an alternative way, via using the relative cup-product argument directly.

\paragraph{Organization. }In \Cref{sec:proofs}, we prove the mixed area separation theorem. In \Cref{sec:proof2}, we establish weak and strict interlacing inequalities and prove Corollary \ref{cor:C}.

\section{Mixed area separation}\label{sec:proofs}
We first prove Theorem \ref{thm:B}. The argument combines the generic results of White and Mazurowski–Zhou and conformal metric perturbation. We begin by fixing the convention and recalling the definition of constrained nondegeneracy.
\begin{convention}[Variation, Jacobi operator]
    We use the variation convention \begin{equation}\label{eq:firstvariation}
 \delta\Area_g(\Sigma)[u]=\int_\Sigma Hu\,d\mu_g,
 \qquad \delta\Vol_g(\Omega)[u]=\int_\Sigma u\,d\mu_g,
\end{equation}
where $u$ is the normal outward variation velocity. We also denote the \emph{Jacobi operator} by: 
\[
 \mathcal L_\Sigma=-\Delta_\Sigma-(|A_\Sigma|^2+\Ric_g(\nu,\nu)).
\]Thus the derivative of scalar mean curvature under a normal variation is \(\mathcal L_\Sigma u\).
\end{convention}

\begin{definition}[Constrained nondegeneracy]\label{def:nd}
\(\Sigma\) is \emph{constrained nondegenerate} if
\[
 \mathcal L_\Sigma u\equiv \text{const.}\quad\text{ on }\Sigma,
 \qquad \int_\Sigma u\,d\mu_g=0
 \quad\Longrightarrow\quad u=0.
\]
The constant is global when \(\Sigma\) is disconnected.
\end{definition}

For reader's convenience, we restate the result of Theorem \ref{thm:B}.

\begin{theorem}\label{thm:nonres}
Let $M^{n+1}$ be a smooth, closed, connected manifold with $n+1\geq 3$. Then there exists a residual set $\mathcal R\subset\Met(M)$ such that, for every $g\in\mathcal R$ and every smooth embedded half-volume CMC boundary $\Sigma=\partial\Omega$,
\begin{equation}\label{eq:rational}
\Area_g(\Sigma)\notin\operatorname{span}_{\mathbb Q}\bigl\{\Area_g(\Gamma):\Gamma\text{ is a smooth closed embedded minimal hypersurface}\bigr\}.
\end{equation}
\end{theorem}

\begin{proof}We first obtain local continuation results for minimal and CMC hypersurfaces under conformal perturbations. Then, we break each fixed rational relation and use a countable covering argument to exclude all such relations simultaneously.

\medskip\noindent
\emph{Step 1. }By White's well known bumpy metrics theorem \cite{WhiteBumpy} and results of Mazurowski--Zhou \cite[Appendix~B]{MZCMC}, there exists a residual set $\mathcal G\subset\Met(M)$ such that 
\begin{enumerate}[label=(\roman*)]
    \item Every smooth closed embedded minimal hypersurface is nondegenerate.
    \item Half-volume minimal hypersurfaces do not exist.
    \item Every smooth almost-embedded half-volume CMC is constrained nondegenerate, embedded, and isolated.
\end{enumerate}
(iii), equivalently, is that the following operator is an isomorphism:
\[
\begin{aligned}
\mathscr A_\Sigma:
C^{2,\alpha}(\Sigma)\times\mathbb R&\longrightarrow C^{0,\alpha}(\Sigma)\times\mathbb R,\\
(u,a)&\longmapsto\left(\mathcal L_\Sigma u-a,\int_\Sigma u\,d\mu_g\right)
\end{aligned}
\]
Indeed, it is a Fredholm operator of index zero and has trivial kernel. Implicit function theorem therefore would give a unique local continuation
\[
g'\longmapsto (\Omega_{g'},\Sigma_{g'},H_{g'})
\]
through half-volume CMC boundaries for all $g'$ sufficiently close to $g$. After shrinking neighborhood of $g$ if necessary, we may always assume $H_{g'}\neq0$. Likewise, same result holds for nondegenerate embedded minimal hypersurfaces, denoted by $g'\mapsto(\Gamma_{g'},H_{g'})$; see \cite[Proposition~B.2 and Corollary~B.4]{MZCMC} and \cite{WhiteSpace} for further discussion.

\medskip\noindent
\emph{Step 2. }Fix finitely many nondegenerate smooth, closed, embedded minimal hypersurfaces $\Gamma_1,\ldots,\Gamma_r$ and $q_1,\ldots,q_r\in\mathbb Q$. By the preceding step, the function
\[
F(g')=\Area_{g'}(\Sigma_{g'})-\sum_{j=1}^r q_j\Area_{g'}(\Gamma_{j,g'})
\]
is well-defined and continuously differentiable on a sufficiently small open neighborhood $U$ of $g$. We claim that $dF\neq 0$ on $U$. It suffices to verify this for an arbitrary metric in $U$ (again, denoted by $g$). Since $\Omega$ is open and the union $\Gamma_1\cup\cdots\cup\Gamma_r$ is finite, we could choose an open ball
\[
B\Subset\Omega\setminus(\Gamma_1\cup\cdots\cup\Gamma_r)
\] and a nonnegative perturbation function $\varphi\in C_c^\infty(B)\setminus \{0\}$. 

Now we can consider the conformal variation $g_t=e^{2t\varphi}g.$ The metric is unchanged near $\Sigma$ and every $\Gamma_j$. In particular, by uniqueness of the local continuation, we have $\Gamma_{j,g_t}=\Gamma_j$ and each $\Gamma_j$ remains minimal for $g_t$. 

Let $\Sigma_t=\partial\Omega_t$ be the corresponding half-volume CMC continuation, and let $u$ denote its outward normal velocity at $t=0$. Since $\varphi$ vanishes near $\Sigma$, we have:
\begin{equation}\label{eq:CMC-variation}
\left.\frac{d}{dt}\right|_{t=0}\Area_{g_t}(\Sigma_t)=H\int_\Sigma u\,d\mu_g=-\frac{n+1}{2}H\int_\Omega\varphi\,dV_g.
\end{equation}
Moreover, differentiating the identity\[
\Vol_{g_t}(\Omega_t)=\frac12\Vol_{g_t}(M)
\] would give \[
\int_\Sigma u\,d\mu_g+(n+1)\int_\Omega\varphi\,dV_g=\frac{n+1}{2}\int_M\varphi\,dV_g.
\]
Since $\operatorname{spt}\varphi\Subset\Omega$, we obtain
\begin{equation}\label{eq:volumevariationmixed}
\int_\Sigma u\,d\mu_g
=-\frac{n+1}{2}\int_\Omega\varphi\,dV_g.
\end{equation}
Combine equations \eqref{eq:CMC-variation}, \eqref{eq:volumevariationmixed} with the vanishing of $H_{\Gamma_t}$ yields
\begin{equation}\label{eq:keyderivative}
\left.\frac{d}{dt}\right|_{t=0}F(g_t)=-\frac{n+1}{2}H\int_\Omega\varphi\,dV_g\neq0.
\end{equation}
The last inequality follows from $H\neq0$, $\varphi\geq0$ and $\varphi\not\equiv0$.

Thus every rational area relation can be broken by a small metric perturbation supported away from all the CMC and minimal hypersurfaces involved.

\medskip\noindent
\emph{Step 3. }It remains to exclude all rational area relations simultaneously. For each $r\geq1$, consider the space of configurations 
\[
\mathcal C_r=
\left\{(g,\Omega,\Sigma,\Gamma_1,\ldots,\Gamma_r):
\begin{array}{l}
g\in\mathcal G,\\
\Sigma=\partial\Omega
\text{ is smooth embedded half-volume CMC boundary},\\
(\Gamma_1,\ldots,\Gamma_r) \text{ is ordered $r$-tuple of smooth embedded}\\ \text{minimal hypersurfaces}.
\end{array}
\right\}.
\]
The space $\mathcal C_r$, equipped with the $C^\infty$ topology on metrics and hypersurfaces and $L^1$ topology on regions, is second countable (hence Lindelöf) for each positive integer $r$; see \cite{DT,WhiteSpace}. At every such configuration, Step~1 would provide a neighborhood in which all the hypersurfaces are uniquely determined by the metric. Consequently, for each $r$, a countable collection of such small neighborhoods covers the full space $\mathcal C_r$. Write $U_i\subset\Met(M)$ for the corresponding metric domains of each charts in $\mathcal C_r$. 

For every rational tuple $q=(q_1,\ldots,q_r)\in\mathbb Q^r$, define on each corresponding chart
\[
F_{i,q}(g')=\Area_{g'}(\Sigma_{g'})-\sum_{j=1}^r q_j\Area_{g'}(\Gamma_{j,g'}).
\]
Its zero set
\[
Z_{i,q}=\{g'\in U_i:F_{i,q}(g')=0\}
\]
is relatively closed in $U_i$. By Step~2, $Z_{i,q}$ has empty interior and is nowhere dense in $U_i$. Since $U_i$ is open in $\Met(M)$, its closure $\overline{Z_{i,q}}$ in $\Met(M)$ is also nowhere dense. Since there are only countably many charts, integers $r\geq1$, and tuples $q\in\mathbb Q^r$. It follows that
\[
\mathcal R:=\mathcal G\setminus
\bigcup_{r=1}^{\infty}\bigcup_i\bigcup_{q\in\mathbb Q^r}\overline{Z_{i,q}}
\]
is residual in $\Met(M)$.

Finally, suppose that $g\in\mathcal R\subset \mathcal G$ and that a smooth embedded half-volume CMC boundary $\Sigma$ satisfies
\[
\Area_g(\Sigma)=\sum_{j=1}^r q_j\Area_g(\Gamma_j)
\]
for some smooth closed embedded minimal hypersurfaces $\Gamma_1,\ldots,\Gamma_r$ and tuple $q_1,\ldots,q_r$. Their combined configuration therefore belongs to one of the continuation charts above, and the identity implies $g\in Z_{i,q}$ for the corresponding configuration chart, a contradiction! Hence \eqref{eq:rational} holds for every $g\in\mathcal R$.
\end{proof}

\section{Interlacing of volume and half-volume spectra}\label{sec:proof2}

We now prove Theorem \ref{thm:A}. The weak inequality follows from a relative cup-product argument and the half-volume correction of \cite{MZCMC}. The strict inequalities then follow from Theorem \ref{thm:B} and the min-max realization theorems from \cite{Law,MZCMC}.

\begin{theorem}\label{thm:mainprecise}
Let $(M^{n+1}, g)$ be a smooth, closed, connected Riemannian manifold with $n+1\geq 3$, then:
\begin{enumerate}[label=(\roman*)]
    \item $\omega_p(M,g)\leq\hw_p(M,g)\leq\omega_{p+1}(M,g)$ for each $p\ge 1$;
    \item Let $3\le n+1\le 5$. There exists a residual set $\mathcal G'\subset\Met(M)$, such that $\omega_p(M,g)<\hw_p(M,g)<\omega_{p+1}(M,g)$ for each $p\ge 1$ and $g\in \mathcal G'$.
\end{enumerate}
\end{theorem}
\begin{proof}
The inequality \(\omega_p\leq\widetilde\omega_p\) follows immediately from the definitions of $\omega_p$ and $\widetilde\omega_p$. To prove the reverse comparison with an index shift, fix \(p\geq1\) and \(\varepsilon>0\). For each integer $p$, let $X$ be a finite cubical subcomplex and let $\Phi :X\longrightarrow\mathcal B(M;\mathcal F;\mathbb Z_2)$ be a near-optimal $p+1$-sweepout, that is: $\sup_{x\in X}\mathbf M(\Phi(x))<\omega_{p+1}+\varepsilon$ for some small $\varepsilon>0$. We shall show that the ordinary sweepout $\Phi$ may naturally induce a half volume $p$-sweepout. 

\medskip\noindent\emph{Step 1. }Let $\alpha:=\Phi^*\lambda\in H^1(X;\mathbb Z_2)$. Since $\Phi$ is a $p+1$-sweepout, by the definition of sweepout and the fact that cup products commutes with pull-back, one has $\alpha^{p+1}\neq 0$. For the rest of proof, we aim to find a subcomplex $P_\delta\subset X$ such that:
\begin{itemize}
    \item  $\alpha^p\big|_{P_\delta}\neq 0$.
    \item $|\Vol(\Omega_x)-\mathfrak h|\le \delta$. Here $\Omega_x$ is either lift to $\mathcal C(M)$ of $\Phi(x)$, i.e. $\partial \Omega_x=\Phi(x)$.
\end{itemize}
For brevity we define the half-volume distance function by:\begin{align*}
    d(T)&:=|\Vol (\Omega)-\mathfrak h|, &&\text{for $T=\partial \Omega\in\mathcal B_n(M;\mathbb Z_2)$.}\\
    d(x)&:=d(\Phi(x)),&&\text{for $x\in X$}.
\end{align*}
We also set $Z=\{x\in X:d(x)=0\}=\Phi^{-1}(\mathcal H_n(M;\mathbb Z_2))\subset X$ be the inverse of zero-level set.

\medskip\noindent\emph{Step 2. }Fix $\delta>0$. Taking any cell and their associated faces of the cubical complex $X$ that intersects $Z$ would give the subcomplex $P_\delta$ with $Z\subset P_\delta\subset X$. Moreover, we may also assume that $d(x)\le\delta$ for all $x\in P_\delta$ by taking possible further finer subdivision of $X$. 

\medskip\noindent\emph{Step 3. }Suppose that $\alpha^p\big|_{P_\delta}= 0$. Pick a small open neighborhood $U\subset X$ such that $U$ deformation retracts to ${P_\delta}$. Thus, we have $H^p(U;\mathbb Z_2)\cong H^p(P_\delta;\mathbb Z_2)$ and $\alpha^p\big|_{U}= 0$ by assumption. 

Take $V:=X\setminus Z$, hence we have:\begin{itemize}
    \item $X=U\cup V$;
    \item $d(x)>0$ for all $x\in V$;
    \item $\Vol(\Omega_x)\neq \mathfrak h$ for all $x\in V$, here $\Omega_x$ is either lifting of $\Phi(x)$ to $\mathcal C(M)$.
\end{itemize}

For each \(x\in V\), let \(\Omega_x^-\) be the unique filling of \(\Phi(x)\) which enclose volume less than \(\mathfrak h\), and define
\[
    s:V\longrightarrow\mathcal C(M),\qquad s(x)=\Omega_x^-.
\]
The double cover $\partial $ admits continuous local lifts. Since local lift agreeing with $\Omega_x^-$ at $x$ still has volume less than $\mathfrak h$ in small neighborhood of $x$, by the uniqueness of the choice of lesser volume filling, the map $s$ is continuous. Write \(\Phi_V:=\Phi|_V\). Then \(\partial\circ s=\Phi_V\), and the following diagram commutes:\[\begin{tikzcd}[row sep=large, column sep=large]
    & \mathcal C(M) \arrow[d, "\partial"] \\
    V \arrow[ur, "s"] 
      \arrow[r, "\Phi_V"'] 
    & \mathcal B(M;\mathcal F;\mathbb Z_2)
\end{tikzcd}\]

Let \(w_1(\partial)\) denote the characteristic class of the double cover $\mathcal C(M)\to\mathcal B(M;\mathcal F;\mathbb Z_2)$. We claim that $ w_1(\partial)=\lambda.$ Indeed, since
\[
    H^1\bigl(\mathcal B(M;\mathcal F;\mathbb Z_2);\mathbb Z_2\bigr)=\{0,\lambda\},
\]
if \(w_1(\partial)=0\), the double cover is trivial and hence $\mathcal C(M)\cong\mathcal B(M;\mathcal F;\mathbb Z_2)\times\mathbb Z_2$. The right-hand side is disconnected while \(\mathcal C(M)\) is contractible, a contradiction.%加引用，可缩

We now consider the pull-back double cover
\[
    E_V:=\bigl\{(x,\Omega)\in V\times\mathcal C(M):\partial\Omega=\Phi_V(x)\bigr\},
    \qquad\pi_V:E_V\longrightarrow V,\quad \pi_V(x,\Omega)=x.
\]
The map \(s\) would give a lift of continuous section
\[
    \widetilde s:V\longrightarrow E_V,
    \qquad
    \widetilde s(x)=(x,s(x)).
\]
Therefore \(\pi_V\) is a trivial double cover, and hence \(w_1(\pi_V)=0\). By naturality of the characteristic class,
\[
    \alpha|_V=\Phi_V^*\lambda=\Phi_V^*w_1(\partial)=w_1(\pi_V)=0
    \quad\text{in }H^1(V;\mathbb Z_2).
\]

\medskip\noindent\emph{Step 4. }We next show that the following items cannot hold simultaneously:
\[
\alpha^{p+1}\neq 0,\qquad\alpha^p\big|_U=0,\qquad\alpha\big|_V=0,\qquad X=U\cup V.
\]
Indeed, if we consider the pair $(X,U)$, the long exact sequence of cohomology contains\[
H^p(X,U;\mathbb Z_2)\xrightarrow{j_U}H^p(X;\mathbb Z_2)\xrightarrow{r_U}H^p(U;\mathbb Z_2)
\]
Since $r_U(\alpha^p)=\alpha^p\big|_U=0$, the exactness would give $\alpha^p\in\operatorname{ker}(r_U)=\operatorname{Im}(j_U)$. Thus there exists $\beta\in H^p(X,U;\mathbb Z_2)$ with $j_U(\beta)=\alpha^p$. Similarly, consider the pair $(X,V)$ and sequence\[
H^1(X,V;\mathbb Z_2)\xrightarrow{j_V}H^1(X;\mathbb Z_2)\xrightarrow{r_V}H^1(V;\mathbb Z_2)
\]would give an element $\gamma\in H^1(X,V;\mathbb Z_2)$ with $j_V(\gamma)=\alpha$.

\medskip\noindent\emph{Step 5. }
Since \(U\) and \(V\) are open subsets of \(X\) and $X=U\cup V$, there is a relative cup product
\[
    H^p(X,U;\mathbb Z_2)\times H^1(X,V;\mathbb Z_2)\xrightarrow{\;\smile\;}H^{p+1}(X,U\cup V;\mathbb Z_2)=H^{p+1}(X,X;\mathbb Z_2)=0.
\]
In particular, \(\beta\smile\gamma=0\).

Let
\[
    j_{U\cup V}:H^{p+1}(X,U\cup V;\mathbb Z_2)\longrightarrow H^{p+1}(X;\mathbb Z_2)
\]
be the canonical map. Compatibility of the relative cup product would give
\[
    0=j_{U\cup V}(\beta\smile\gamma)=j_U(\beta)\smile j_V(\gamma)=\alpha^p\smile\alpha=\alpha^{p+1},
\]
a contradiction! Therefore $\alpha^p|_{P_\delta}\neq0.$

\medskip\noindent\emph{Step 6. }Define
\[
    \Phi_\delta:=\Phi\big|_{P_\delta}:P_\delta\longrightarrow\mathcal B(M;\mathcal F;\mathbb Z_2).
\]
By naturality of pull-back and restriction,
\[
    \Phi_\delta^*(\lambda^p)=\bigl(\Phi^*(\lambda^p)\bigr)\big|_{P_\delta}=\alpha^p\big|_{P_\delta}
    \neq0.
\]
This map is naturally $\mathcal F$-continuous and has no concentration of mass since it is a restriction of a $p+1$-sweepout $\Phi$. Hence \(\Phi_\delta\) is an ordinary \(p\)-sweepout. Moreover, the construction of \(P_\delta\) would give
\[
    d(\Phi_\delta(x))\leq\delta\qquad\text{for every }x\in P_\delta,
\]
and
\[
    \sup_{x\in P_\delta}\mathbf M(\Phi_\delta(x))\leq\sup_{x\in X}\mathbf M(\Phi(x))<\omega_{p+1}+\varepsilon.
\]

We now try to reduce this near half-volume $p$-sweepout to an exact half-volume $p$-sweepout, the main tool is the following volume correction retraction given by Mazurowski-Zhou:%加引用

\begin{proposition}[Cf. {\cite{half} and \cite[Lemma 6.6]{MZCMC}}]\label{prop:external}%引用
    There exists a deformation retraction $\theta:\mathcal B(M;\mathbb Z_2)\times [0,1]\longrightarrow\mathcal B(M;\mathbb Z_2)$ such that \begin{enumerate}[label=(\roman*)]
    \item $\theta(T,1)\in\mathcal H(M;\mathbb Z_2)$;
    \item $\theta(T,t)=\partial (\Omega\cup A_{s(\Omega,t)})$, where $\Omega$ is the filling of $T$ with $\Vol(\Omega)\leq\mathfrak h$, and \(A_s\) is the sub-level set which enclose volume $s$ of some Morse function $f:M\rightarrow \mathbb R$.
    \item $\mathbf M(\theta(T,1))\leq\mathbf M(T)+w(d(T))$ for all $T\in\mathcal B(M;\mathbb Z_2).$ Here $w:[0,+\infty)\to [0,+\infty)$ is a continuous function with $w(0)=0$.
    \end{enumerate}
\end{proposition}
\begin{proof}
    The map $\theta$ was first constructed in \cite{half}, while the function $w$ was given by \cite[Lemma 6.6]{MZCMC}.
\end{proof}
Given above deformation retraction, one may define the half-volume $p$-sweepout by \[
\Psi_\delta:P_\delta\longrightarrow\mathcal H(M;\mathbb Z_2),\qquad x\longmapsto\theta(\Phi_\delta(x),1).
\]Indeed, since $\Phi_\delta$ is a $p$-sweepout and $\theta $ is a deformation retraction, we get ${\Psi_\delta}^*\lambda^p\neq 0$. $\Psi_\delta$ has no concentration of mass since \[
\begin{aligned}
\sup_{x,q}\|\Psi_\delta(x)\|(B_r(q))\le {}&
\sup_{x,q}\|\Phi_\delta(x)\|(B_r(q))\\
&+\sup_{s,q}\mathcal H^n(\partial^*A_s\cap B_r(q)).
\end{aligned}
\]The right hand side has no-concentration property uniformly. Proposition \ref{prop:external} therefore give $\tilde{\omega}_p\leq\omega_{p+1}+\varepsilon+\sup_{0\leq t\leq\delta}w(t)$. Let \(\delta\downarrow0\) and then \(\eps\downarrow0\), this proves the first assertion.

\medskip\noindent\emph{Step 7. }For dimensions $3\le n+1\le 5$ and a residual set $\mathcal G$ of smooth metrics, each $\omega_p$ is realized by area of a finite integral sum of areas of smooth closed embedded minimal hypersurfaces, and each $\tilde{\omega_p}$ is realized by a smooth, embedded CMC hypersurface; see \cite{Law,MZCMC}. Let $\mathcal G'=\mathcal G\cap\mathcal R$, where $\mathcal R$ is the residual set from Theorem \ref{thm:nonres}, the second assertion follows immediately.\end{proof}

Mazurowski--Zhou \cite{MZCMC} proved the infinite existence of half-volume CMC hypersurfaces for generic metrics. Theorem \ref{thm:mainprecise} places a half-volume CMC boundary in each interval $(\omega_p,\omega_{p+1})$. Combining this with the Weyl law gives the quantitative lower bound, as in the following corollary.

\begin{corollary}\label{cor:Cprecise}
For every $g\in\mathcal G'$, let $\mathcal N_g(A)$ denote the number of distinct smooth embedded half-volume CMC boundaries with nonzero mean curvature and area at most $A$, then
\[
\liminf_{A\to\infty}
\frac{\mathcal N_g(A)}{A^{n+1}}\ge
\frac{1}{a_n^{n+1}\Vol_g(M)^n},
\]
where $a_n$ is the dimensional constant in \cite[Theorem 1.1]{Law}.
\end{corollary}

\begin{proof}
By Theorem~\ref{thm:mainprecise}, especially step 7, the hypersurfaces $\Sigma_p$ are pairwise distinct and satisfy
\[
\omega_p(M,g)< \Area_g(\Sigma_p)< \omega_{p+1}(M,g).
\]
Consequently,
\[
\mathcal N_g(A)\ge
\max\bigl\{0,\#\{p\ge1:\omega_p(M,g)\le A\}-1\bigr\}.
\]
The conclusion follows immediately from \cite[Theorem~1.1]{Law}.
\end{proof}

\begin{flushleft}
\medskip\noindent
\begin{tabbing}

				Xiaoxiang Jiao\\
				School of Mathematical Sciences, University of Chinese Academy of Sciences\\
				19A Yuquan Road, Beijing, 100049, China\\
				\texttt{xxjiao@ucas.ac.cn}
				
\end{tabbing}

\begin{tabbing}
Hangyue Zhu\\School of Mathematical Sciences, University of Chinese Academy of Sciences\\
				19A Yuquan Road, Beijing, 100049, China\\
				\texttt{zhuhangyue24@mails.ucas.ac.cn}
\end{tabbing}

\end{flushleft}

\paragraph{Conflict of Interests.}
On behalf of all authors, the corresponding author states that there is no conflict of interest.

\paragraph{Data Availability.}
Data sharing is not applicable to this article, as no datasets were generated or analyzed during the current study.

\paragraph{Funding.}
The authors are supported by the National Natural Science Foundation of China (Grant no. 12371055). 

\paragraph{Acknowledgments.}
The present work was motivated in part by \cite[Remark 25]{half}, which also suggested the proof strategy of Theorem \ref{thm:mainprecise} (i).

\paragraph{Declaration of Generative AI Use.} During the preparation of this work, the authors used ChatGPT and Deepseek to improve the grammar, clarity, and quality of the manuscript's writing, and to assist in creating TikZ figure. After using these tools, the authors reviewed and edited the content as needed and take full responsibility for the content of the article.%这里单复数下次文章复制记得改

\bibliographystyle{alpha}
\bibliography{ref}

@misc{Talipov,
      title={Equivalence of Almgren-Pitts and phase-transition half-volume spectra}, 
      author={Talant Talipov},
      year={2026},
      eprint={2604.01091},
      archivePrefix={arXiv},
      primaryClass={math.DG},
        howpublished = {\href{https://arxiv.org/abs/2604.01091}{arXiv:2604.01091}},
}

@article {CMCI,
    AUTHOR = {Zhou, Xin and Zhu, Jonathan J.},
     TITLE = {Min-max theory for constant mean curvature hypersurfaces},
   JOURNAL = {Invent. Math.},
  FJOURNAL = {Inventiones Mathematicae},
    VOLUME = {218},
      YEAR = {2019},
    NUMBER = {2},
     PAGES = {441--490},
      ISSN = {0020-9910,1432-1297},
   MRCLASS = {53C42 (28A75 49J35 49Q15 53A15)},
  MRNUMBER = {4011704},
MRREVIEWER = {John\ McCuan},
       DOI = {10.1007/s00222-019-00886-1},
       URL = {https://doi.org/10.1007/s00222-019-00886-1},
}

@article {MR2481045,
    AUTHOR = {Pacard, F. and Xu, X.},
     TITLE = {Constant mean curvature spheres in {R}iemannian manifolds},
   JOURNAL = {Manuscripta Math.},
  FJOURNAL = {Manuscripta Mathematica},
    VOLUME = {128},
      YEAR = {2009},
    NUMBER = {3},
     PAGES = {275--295},
      ISSN = {0025-2611,1432-1785},
   MRCLASS = {53C42},
  MRNUMBER = {2481045},
MRREVIEWER = {Alma\ Luisa\ Albujer Brotons},
       DOI = {10.1007/s00229-008-0230-7},
       URL = {https://doi.org/10.1007/s00229-008-0230-7},
}

@article {Dey222,
    AUTHOR = {Dey, Akashdeep},
     TITLE = {A comparison of the {A}lmgren-{P}itts and the {A}llen-{C}ahn
              min-max theory},
   JOURNAL = {Geom. Funct. Anal.},
  FJOURNAL = {Geometric and Functional Analysis},
    VOLUME = {32},
      YEAR = {2022},
    NUMBER = {5},
     PAGES = {980--1040},
      ISSN = {1016-443X,1420-8970},
   MRCLASS = {53A10 (58E12)},
  MRNUMBER = {4498838},
MRREVIEWER = {Erasmo\ Caponio},
       DOI = {10.1007/s00039-022-00610-x},
       URL = {https://doi.org/10.1007/s00039-022-00610-x},
}

@misc{R6,
      title={Do Carmo's problem for CMC hypersurfaces in $\mathbb{R}^6$}, 
      author={Jingche Chen and Han Hong and Haizhong Li},
      year={2025},
      eprint={2503.08107},
      archivePrefix={arXiv},
      primaryClass={math.DG},
      howpublished = {\href{https://arxiv.org/abs/2503.08107}{arXiv:2503.08107}},
}

@misc{MZPisa,
      title={An Alternative for Constant Mean Curvature Hypersurfaces}, 
      author={Liam Mazurowski and Xin Zhou},
      year={2024},
      eprint={2408.13864},
      archivePrefix={arXiv},
      primaryClass={math.DG},
      howpublished = {\href{https://arxiv.org/abs/2408.13864}{arXiv:2408.13864}},
      note={To appear at \emph{Ann. Scuola Norm. Sup. Pisa Cl. Sci.}}
}

@misc{MZCMC,
      title={Infinitely Many Half-Volume Constant Mean Curvature Hypersurfaces via Min-Max Theory}, 
      author={Liam Mazurowski and Xin Zhou},
      year={2024},
      eprint={2405.00595},
      archivePrefix={arXiv},
      primaryClass={math.DG},
      howpublished = {\href{https://arxiv.org/abs/2405.00595v2}{arXiv:2405.00595v2}},
      note={To appear at \emph{J. Differential Geom.}}
}

@article {half,
    AUTHOR = {Mazurowski, Liam and Zhou, Xin},
     TITLE = {The half-volume spectrum of a manifold},
   JOURNAL = {Calc. Var. Partial Differential Equations},
  FJOURNAL = {Calculus of Variations and Partial Differential Equations},
    VOLUME = {64},
      YEAR = {2025},
    NUMBER = {5},
     PAGES = {Paper No. 155, 22},
      ISSN = {0944-2669,1432-0835},
   MRCLASS = {53A10},
  MRNUMBER = {4902143},
MRREVIEWER = {Nicolau\ S.\ Aiex},
       DOI = {10.1007/s00526-025-02949-z},
       URL = {https://doi.org/10.1007/s00526-025-02949-z},
}

@article{ALMGREN1962257,
title = {The homotopy groups of the integral cycle groups},
journal = {Topology},
volume = {1},
number = {4},
pages = {257-299},
year = {1962},
issn = {0040-9383},
doi = {https://doi.org/10.1016/0040-9383(62)90016-2},
url = {https://www.sciencedirect.com/science/article/pii/0040938362900162},
author = {Frederick Justin Almgren}
}

@article {DeyCMC,
    AUTHOR = {Dey, Akashdeep},
     TITLE = {Existence of multiple closed {CMC} hypersurfaces with small
              mean curvature},
   JOURNAL = {J. Differential Geom.},
  FJOURNAL = {Journal of Differential Geometry},
    VOLUME = {125},
      YEAR = {2023},
    NUMBER = {2},
     PAGES = {379--403},
      ISSN = {0022-040X,1945-743X},
   MRCLASS = {53C42 (53A10)},
  MRNUMBER = {4649390},
MRREVIEWER = {Fei-Tsen\ Liang},
       DOI = {10.4310/jdg/1696432925},
       URL = {https://doi.org/10.4310/jdg/1696432925},
}

@book{DT,
  title={Differential Topology},
  author={Hirsch, M.W.},
  isbn={9781468494495},
  series={Graduate Texts in Mathematics},
  url={https://books.google.com/books?id=emTmBwAAQBAJ},
  year={2012},
  publisher={Springer New York}
}

@article {WhiteSpace,
    AUTHOR = {White, Brian},
     TITLE = {The space of minimal submanifolds for varying {R}iemannian
              metrics},
   JOURNAL = {Indiana Univ. Math. J.},
  FJOURNAL = {Indiana University Mathematics Journal},
    VOLUME = {40},
      YEAR = {1991},
    NUMBER = {1},
     PAGES = {161--200},
      ISSN = {0022-2518,1943-5258},
   MRCLASS = {58D10 (53C42)},
  MRNUMBER = {1101226},
MRREVIEWER = {Jo\~ao\ Lucas Marques Barbosa},
       DOI = {10.1512/iumj.1991.40.40008},
       URL = {https://doi.org/10.1512/iumj.1991.40.40008},
}

@article {WhiteBumpy,
    AUTHOR = {White, Brian},
     TITLE = {On the bumpy metrics theorem for minimal submanifolds},
   JOURNAL = {Amer. J. Math.},
  FJOURNAL = {American Journal of Mathematics},
    VOLUME = {139},
      YEAR = {2017},
    NUMBER = {4},
     PAGES = {1149--1155},
      ISSN = {0002-9327,1080-6377},
   MRCLASS = {53C42},
  MRNUMBER = {3689325},
MRREVIEWER = {Otis\ Chodosh},
       DOI = {10.1353/ajm.2017.0029},
       URL = {https://doi.org/10.1353/ajm.2017.0029},
}

@article {Law,
    AUTHOR = {Liokumovich, Yevgeny and Marques, Fernando C. and Neves,
              Andr\'e},
     TITLE = {Weyl law for the volume spectrum},
   JOURNAL = {Ann. of Math. (2)},
  FJOURNAL = {Annals of Mathematics. Second Series},
    VOLUME = {187},
      YEAR = {2018},
    NUMBER = {3},
     PAGES = {933--961},
      ISSN = {0003-486X,1939-8980},
   MRCLASS = {53C23 (58E05 58J50)},
  MRNUMBER = {3779961},
MRREVIEWER = {Leonid\ Friedlander},
       DOI = {10.4007/annals.2018.187.3.7},
       URL = {https://doi.org/10.4007/annals.2018.187.3.7},
}
\end{document}